\documentclass[eqno]{amsart}
\usepackage[dvips,final]{epsfig}
\usepackage{amsmath}
\usepackage{amscd}
\usepackage{amsfonts}
\usepackage{amssymb}
\usepackage{latexsym}
\usepackage{color}
\title{Bi-Lipschitz $\mathcal{A}$-equivalence of $\mathcal{K}$-equivalent map-germs}
\author{
J.~C.~F.~Costa}
\address{
Departamento de Matem\'atica, IBILCE-UNESP, Campus de S\~ao Jos\'e
do Rio Preto SP, BRAZIL. } \email{jcosta@ibilce.unesp.br}
\author{T.~Nishimura}
\address{
Research Group of Mathematical Sciences, Research Institute of
Environment and Information Sciences,
Yokohama National University,
Yokohama 240-8501, JAPAN. } \email{nishimura-takashi-yx@ynu.jp}
\thanks{T.~Nishimura was partially supported
by JSPS and CAPES under the Japan--Brazil research cooperative
program}
\author{M.~A.~S.~Ruas}
\address{
Departamento de Matem\'atica, Instituto de Ci\'encias Matem'aticas e
de Computac\~ao, Universidade de S\~ao Paulo, Caixa Postal 668,
13560-970 S\~ao Carlos SP, BRAZIL. } \email{maasruas@icmc.sc.usp.br}
\thanks{The authors J.~C.~F. Costa and M.~A.~S.~Ruas  were partially supported by CNPq, CAPES and
FAPESP}
\dedicatory{Dedicated to Professor Heisuke Hironaka on the occasion of his 80th birthday}
 \keywords{bi-Lipschitz
$\mathcal{A}$-equivalence, Lipschitz trivial, Lipschitz version of
Fukuda-Fukuda theorem. }
\begin{document}
\begin{abstract}\quad
In this paper, two sufficient conditions are provided for given two
$\mathcal{K}$-equivalent map-germs to be bi-Lipschitz
$\mathcal{A}$-equivalent. These are Lipschitz analogues of the known
results on $C^r$ $\mathcal{A}$-equivalence $(0\le r\le \infty)$ for
given two $\mathcal{K}$-equivalent map-germs. As a corollary of one
of our results, a Lipschitz version of the well-known Fukuda-Fukuda
theorem is provided.
\end{abstract}
\date{}
\maketitle \noindent
\newtheorem{theorem}{Theorem}
\newtheorem{corollary}{Corollary}
\newtheorem{lemma}{Lemma}[section]
\newtheorem{proposition}{Proposition}
\newtheorem{assertion}{Assertion}[section]
\newtheorem{definition}{Definition}
\newtheorem{example}{Example}[section]
\newtheorem{conjecture}{Conjecture}
\def\theenumi{\roman{enumi}}
\noindent
\section{Introduction} 
\label{intro} \noindent The classification problem of generic
singularities of $C^\infty$ mappings up to $C^\infty$ coordinate
transformations is one of the most important problems in Singularity
Theory.    However, since it is impossible to obtain countably many
list on this problem, the classification of singularities of
$C^\infty$ mappings by weaker equivalence relations are often
studied. The $C^r$ $\mathcal{A}$-equivalence $(0\le r < \infty)$ has
been well-studied. { } The bi-Lipschitz
$\mathcal{A}$-equivalence is also one of the weaker equivalence
relations and the classification problem with respect to
bi-Lipschitz $\mathcal{A}$-equivalence is an interesting subject in
the recent development of metric singularity theory (see for
instance \cite{costasaiasoares, fernandesruas, fernandessoares,
henryparusinski}).
In this paper we give sufficient conditions for the bi-Lipschitz
$\mathcal{A}$-equivalence of $\mathcal{K}$-equivalent map-germs,
providing Lipschitz versions of results with respect to the $C^r$
$\mathcal{A}$-equivalence given in \cite{nishimuralondon,
nishimuratopology} where $0\le r\le \infty$. { }
\par
A map $\varphi: U\subset \mathbb{R}^n\to \mathbb{R}^p$ is said to be
{\it Lipschitz} if there exists a constant $c > 0$ such that the
following holds, where $U$ is an open set of $\mathbb{R}^n$.
\[
||\varphi(x)-\varphi(y)||\le c ||x-y|| \quad \mbox{for any }x,y\in
U.
\]
A Lipschitz map $\varphi : U\subset \mathbb{R}^n\to
\varphi(U)\subset\mathbb{R}^p$ is called a {\it bi-Lipschitz
homeomorphism} if $n=p$ and $\varphi$ has a Lipschitz inverse. Two
$C^\infty$ map-germs $f, g:(\mathbb{R}^n,0)\to (\mathbb{R}^p,0)$ are
said to be {\it bi-Lipschitz} $\mathcal{A}$-{\it equivalent} if
there exist germs of bi-Lipschitz homeomorphisms $s:
(\mathbb{R}^n,0)\to (\mathbb{R}^n,0)$ and $t: (\mathbb{R}^p,0)\to
(\mathbb{R}^p,0)$ such that $f=t\circ g\circ s$.
\par
Let $f: (\mathbb{R}^n,0)\to (\mathbb{R}^p,0)$ be a $C^\infty$
map-germ. Any $C^\infty$ map-germ $\Phi: (\mathbb{R}^n\times
\mathbb{R}^k, (0,0))\to (\mathbb{R}^p,0)$ such that $\Phi(x,0)=f(x)$
is called a {\it $C^\infty$ deformation-germ of $f$}.
\begin{definition}\label{def1}
{\rm Let $f: (\mathbb{R}^n,0)\to (\mathbb{R}^p,0)$ be a $C^\infty$
map-germ. A $C^\infty$ deformation-germ $\Phi: (\mathbb{R}^n\times
\mathbb{R}^k, (0,0))\to (\mathbb{R}^p,0)$ of $f$ is said to be {\it
Lipschitz trivial} if there exist $L$-stratifications in the sense
of \cite{mostowski, parusinski, parusinski2}, $\mathcal{S}$ of
$\mathbb{R}^n\times \mathbb{R}^k$, $\mathcal{T}$ of
$\mathbb{R}^p\times \mathbb{R}^k$ and $\{\mathbb{R}^k\}$ of
$\mathbb{R}^k$ such that the following three conditions are
satisfied:
\begin{enumerate}
\item The map-germ $(\Phi, \pi): (\mathbb{R}^n\times \mathbb{R}^k,(0,0))\to (\mathbb{R}^p\times \mathbb{R}^k,(0,0))$ is a stratified map
with respect to $\mathcal{S}$ and $\mathcal{T}$, where $\pi:
(\mathbb{R}^n\times \mathbb{R}^k,(0,0))\to (\mathbb{R}^k,0)$ is the
canonical projection.
\item The canonical projection $\pi': (\mathbb{R}^p\times \mathbb{R}^k,(0,0))\to (\mathbb{R}^k,0)$ is a stratified map
with respect to $\mathcal{T}$ and $\{\mathbb{R}^k\}$.
\item There exist germs of stratified Lipschitz vector fields $\xi_i+\frac{\partial}{\partial \lambda_i}$ at the origin in
$\mathbb{R}^n\times \mathbb{R}^k$ with respect to $\mathcal{S}$ and
$\eta_i+\frac{\partial}{\partial \lambda_i}$ at the origin in
$\mathbb{R}^p\times \mathbb{R}^k$ with respect to $\mathcal{T}$ such
that
\[
\xi_i+ \frac{\partial}{\partial \lambda_i} \mbox{ lifts } \eta_i+
\frac{\partial}{\partial \lambda_i} \mbox{ with respect to } (\Phi,
\pi),
\]
where $\frac{\partial}{\partial \lambda_i}$ is the trivial vector
field of $\mathbb{R}^k$ with respect to the $i$-th coordinate
function $\lambda_i$ of the standard coordinate neighborhood
$(\mathbb{R}^k, (\lambda_1, \ldots, \lambda_k))$.
\end{enumerate}
}
\end{definition}

{Definition \ref{def1} was motivated by  the definition of {\it Thom
trivial} deformation  given by Nishimura in \cite{nishimuralondon}.
A deformation is Thom trivial   if there exists a $C$-regular
stratification in the sense of K. Bekka (\cite{bekka}), satisfying
conditions (i) and (ii) above where $(\Phi, \pi)$ is a Thom map.
For any Thom trivial deformation we can apply  
Thom's second
isotopy lemma (\cite{bekka, gibsonetall, matherhowto,
duplessiswall})  which implies that $(\Phi, \pi)$ is topologically
equivalent to $(f,\pi)$. To obtain a similar condition to that for
the Lispchitz case, we replace $C$-regular stratifications by
$L$-stratifications (in the sense of Mostowski \cite{mostowski}) and
we add the condition (iii) above, because there seems to have been
no bi-Lipschitz versions of Thom's second isotopy lemma.}
\smallskip
\par
It follows from Definition \ref{def1} that
for any Lipschitz trivial deformation-germ $\Phi:
(\mathbb{R}^n\times \mathbb{R}^k,(0,0))\to (\mathbb{R}^p,0)$ of $f$,
there exist germs of bi-Lipschitz homeomorphisms $h:
(\mathbb{R}^n\times \mathbb{R}^k, (0,0))\to (\mathbb{R}^n\times
\mathbb{R}^k,(0,0))$ and $H: (\mathbb{R}^p\times \mathbb{R}^k,
(0,0))\to (\mathbb{R}^p\times \mathbb{R}^k,(0,0))$ such that the
following diagram is commutative:
\begin{equation}
\begin{CD}
{ } @. (\mathbb{R}^n\times \mathbb{R}^k,(0,0)) @>(\Phi, \pi)>> (\mathbb{R}^p\times\mathbb{R}^k,(0,0)) @>\pi'>>(\mathbb{R}^k,0) \\
@. @VhVV @VHVV @| \\
{ } @. (\mathbb{R}^n\times \mathbb{R}^k,(0,0)) @>(f, \pi)>>
(\mathbb{R}^p\times\mathbb{R}^k,(0,0)) @>\pi'>>(\mathbb{R}^k,0),
\end{CD}
\end{equation}
where $\pi: (\mathbb{R}^n\times \mathbb{R}^k,(0,0))\to
(\mathbb{R}^k,0)$ and $\pi': (\mathbb{R}^p\times
\mathbb{R}^k,(0,0))\to (\mathbb{R}^k,0)$ are canonical projections.
The main result of this paper is the following:
\begin{theorem}\label{thm1}
Let $f, g: (\mathbb{R}^n,0)\to (\mathbb{R}^p,0)$ be two $C^\infty$
map-germs. Suppose that the following three conditions are
satisfied.
\begin{enumerate}
\item There exist a germ of $C^\infty$ diffeomorphism $s: (\mathbb{R}^n,0)\to (\mathbb{R}^n,0)$ and
a $C^\infty$ map-germ $M: (\mathbb{R}^n, 0)\to (GL(p, \mathbb{R}),
M(0))$ such that $f(x)=M(x)g(s(x))$.
\item The $C^\infty$ map-germ $F: (\mathbb{R}^n\times \mathbb{R}^p, (0,0))\to (\mathbb{R}^p,0)$ given by
\[
F(x, \lambda)=f(x)-M(x)\lambda
\]
is a Lipschitz trivial deformation-germ of $f$.
\item The $C^\infty$ map-germ $G: (\mathbb{R}^n\times \mathbb{R}^p, (0,0))\to (\mathbb{R}^p,0)$ given by
\[
G(x, \lambda)=g(x)-M(s^{-1}(x))^{-1}\lambda
\]
is a Lipschitz trivial deformation-germ of $g$.
\end{enumerate}
Then, $f$ and $g$ are bi-Lipschitz $\mathcal{A}$-equivalent.
\end{theorem}

{ } The condition (i) in Theorem \ref{thm1} is equivalent to
say that $f$ and $g$ are $\mathcal{K}$-equivalent (see
\cite{mather3}).{ }
\begin{definition}
{\rm A $C^\infty$ map-germ $f: (\mathbb{R}^n,0)\to (\mathbb{R}^p,0)$
is said to be {\it Lipschitz stable} if any $C^\infty$
deformation-germ of $f$ is Lipschitz trivial. }
\end{definition}
As a corollary of Theorem 1, the following Lipschitz version of
Fukuda-Fukuda theorem \cite{fukudafukuda} is obtained.
\begin{corollary}
Let $f, g: (\mathbb{R}^n,0)\to (\mathbb{R}^p,0)$ be two Lipschitz
stable map-germs. Suppose that there exist a germ of $C^\infty$
diffeomorphism $s: (\mathbb{R}^n,0)\to (\mathbb{R}^n,0)$ and a
$C^\infty$ map-germ $M: (\mathbb{R}^n, 0)\to (GL(p, \mathbb{R}),
M(0))$ such that $f(x)=M(x)g(s(x))$. Then, $f$ and $g$ are
bi-Lipschitz $\mathcal{A}$-equivalent.
\end{corollary}
In the case that both the given two map-germs $f$ and $g$ are of
rank zero, Theorem 1 can be strengthened as follows.
\begin{theorem}\label{thm2}
Let $f, g: (\mathbb{R}^n,0)\to (\mathbb{R}^p,0)$ be two $C^\infty$
map-germs of rank zero. Suppose that the following two conditions
are satisfied.
\begin{enumerate}
\item There exist a germ of $C^\infty$ diffeomorphism $s: (\mathbb{R}^n,0)\to (\mathbb{R}^n,0)$ and
a $C^\infty$ map-germ $M: (\mathbb{R}^n, 0)\to (GL(p, \mathbb{R}),
M(0))$ such that $f(x)=M(x)g(s(x))$.
\item The $C^\infty$ map-germ $F: (\mathbb{R}^n\times \mathbb{R}^p, (0,0))\to (\mathbb{R}^p,0)$ given by
\[
F(x, \lambda)=f(x)-M(x)\lambda
\]
is a Lipschitz trivial deformation-germ of $f$.
\end{enumerate}
Then, $f$ and $g$ are bi-Lipschitz $\mathcal{A}$-equivalent.
\end{theorem}
\bigskip
\par

{ }\emph{\textbf{Remark 1.}} Let $f,g:(\mathbb{R}^2,0) \to
(\mathbb{R}^2,0)$ given by
$$
f(x,y) = (x,y^3+xy) \,\,\,\,\,\,\,{\rm and} \,\,\,\,\,\,\, g(x,y) =
(x,y^3).
$$
It is easily seen that $f(x,y)=M(x,y)g(x,y)$ with
$M:(\mathbb{R}^2,0) \to (GL(2,\mathbb{R}),E_2)$. Moreover, it is
well known that $f$ is infinitesimally stable, hence any $C^\infty$ deformation germ of $f$ is $C^\infty$
trivial (and hence Lipschitz trivial). However, $f$ and $g$ are not
bi-Lipschitz $\mathcal{A}$-equivalent. This example shows that
rank zero hypothesis is essential in Theorem \ref{thm2}.
{ }

\vspace{.5cm}

In Section 2, $C^r$ versions of Theorem 1, Corollary 1 and Theorem 2
are reviewed where $0\le r \le \infty$. In Section 3, a strategy to
show Theorems 1 and 2 is explained. Theorems 1 and 2 are proved in
Section 4.
\section{The $C^r$ versions $(0\le r\le \infty)$}
\label{section 2}
In this section, for the readers' convenience, we review $C^r$
versions of Theorem 1, Corollary 1 and Theorem 2 $(0\le r\le
\infty)$ briefly. For details on $C^r$ versions $(1\le r\le
\infty)$, refer to \cite{nishimuratopology}; and on $C^0$ versions,
refer to \cite{nishimuralondon} (see also \cite{nishimurabanach}
where \cite{nishimuralondon, nishimuratopology} and related results
are summarized).
\begin{definition}
{\rm
\begin{enumerate}
\item Let $f, g:(\mathbb{R}^n,0)\to (\mathbb{R}^p,0)$ be two $C^\infty$ map-germs.
\begin{enumerate}
\item The given $f$ and $g$ are said to be
$C^r$ $\mathcal{A}$-{\it equivalent} $(1\le r\le \infty)$ if there
exist germs of $C^r$ diffeomorphisms $s: (\mathbb{R}^n, 0)\to
(\mathbb{R}^n, 0)$ and $t: (\mathbb{R}^p, 0)\to (\mathbb{R}^p, 0)$
such that $f=t\circ g\circ s$.
\item The given $f$ and $g$ are said to be
{\it topologically} $\mathcal{A}$-{\it equivalent} if there exist
germs of homeomorphisms $s: (\mathbb{R}^n,0)\to (\mathbb{R}^n,0)$
and $t: (\mathbb{R}^p,0)\to (\mathbb{R}^p,0)$ such that $f=t\circ
g\circ s$.
\end{enumerate}
\item Let $f: (\mathbb{R}^n,0)\to (\mathbb{R}^p,0)$ be a $C^\infty$ map-germ.
A $C^\infty$ deformation-germ $\Phi: (\mathbb{R}^n\times
\mathbb{R}^k, (0,0))\to (\mathbb{R}^p,0)$ of $f$ is said to be {\it
$C^r$ $\mathcal{A}$-trivial} $(1\le r\le \infty)$ if there exist
germs of $C^r$ diffeomorphisms $h: (\mathbb{R}^n\times \mathbb{R}^k,
(0,0))\to (\mathbb{R}^n\times \mathbb{R}^k,(0,0))$ and $H:
(\mathbb{R}^p\times \mathbb{R}^k, (0,0))\to (\mathbb{R}^p\times
\mathbb{R}^k,(0,0))$ such that the diagram (1) in Section 1 is
commutative.
\end{enumerate}
}
\end{definition}

{ }
\subsection{The $C^r$ versions $(1\le r \le \infty)$}
\begin{theorem}[\cite{nishimuratopology}]
Let $f, g: (\mathbb{R}^n,0)\to (\mathbb{R}^p,0)$ be two $C^\infty$
map-germs. Suppose that the following three conditions are
satisfied.
\begin{enumerate}
\item There exist a germ of $C^\infty$ diffeomorphism $s: (\mathbb{R}^n,0)\to (\mathbb{R}^n,0)$ and
a $C^\infty$ map-germ $M: (\mathbb{R}^n, 0)\to (GL(p, \mathbb{R}),
M(0))$ such that $f(x)=M(x)g(s(x))$.
\item The $C^\infty$ map-germ $F: (\mathbb{R}^n\times \mathbb{R}^p, (0,0))\to (\mathbb{R}^p,0)$ given by
\[
F(x, \lambda)=f(x)-M(x)\lambda
\]
is a $C^r$ $\mathcal{A}$-trivial deformation-germ of $f$.
\item The $C^\infty$ map-germ $G: (\mathbb{R}^n\times \mathbb{R}^p, (0,0))\to (\mathbb{R}^p,0)$ given by
\[
G(x, \lambda)=g(x)-M(s^{-1}(x))^{-1}\lambda
\]
is a $C^r$ $\mathcal{A}$-trivial deformation-germ of $g$.
\end{enumerate}
Then, $f$ and $g$ are $C^r$ $\mathcal{A}$-equivalent $(1\le r \le
\infty)$.
\end{theorem}
A $C^\infty$ map-germ $f: (\mathbb{R}^n,0)\to (\mathbb{R}^p,0)$ is
said to be $C^r$ {\it stable} $(1\le r \le \infty)$ if any
$C^\infty$ deformation of $f$ is $C^r$ $\mathcal{A}$-trivial $(1\le
r \le \infty)$. It is common in the literature to call a $C^\infty$
stable map-germ just by stable map-germ.

\begin{corollary}[\cite{nishimuratopology}]\label{corol1}
Let $f, g: (\mathbb{R}^n,0)\to (\mathbb{R}^p,0)$ be two $C^r$ stable
map-germs $(1\le r \le \infty)$. Suppose that there exist a germ of
$C^\infty$ diffeomorphism $s: (\mathbb{R}^n,0)\to (\mathbb{R}^n,0)$
and a $C^\infty$ map-germ $M: (\mathbb{R}^n, 0)\to (GL(p,
\mathbb{R}), M(0))$ such that $f(x)=M(x)g(s(x))$. Then, $f$ and $g$
are $C^r$ $\mathcal{A}$-equivalent.
\end{corollary}

Notice that for $r=\infty$,  Corollary \ref{corol1} is the
classical Mather's classification theorem for stable map-germs (see
\cite{mather4}). { }
\begin{theorem}[\cite{nishimuratopology}]
Let $f, g: (\mathbb{R}^n,0)\to (\mathbb{R}^p,0)$ be two $C^\infty$
map-germs of rank zero. Suppose that the following two conditions
are satisfied.
\begin{enumerate}
\item There exist a germ of $C^\infty$ diffeomorphism $s: (\mathbb{R}^n,0)\to (\mathbb{R}^n,0)$ and
a $C^\infty$ map-germ $M: (\mathbb{R}^n, 0)\to (GL(p, \mathbb{R}),
M(0))$ such that $f(x)=M(x)g(s(x))$.
\item The $C^\infty$ map-germ $F: (\mathbb{R}^n\times \mathbb{R}^p, (0,0))\to (\mathbb{R}^p,0)$ given by
\[
F(x, \lambda)=f(x)-M(x)\lambda
\]
is a $C^r$ $\mathcal{A}$-trivial deformation-germ of $f$.
\end{enumerate}
Then, $f$ and $g$ are $C^r$ $\mathcal{A}$-equivalent { }($1
\le r \le \infty$).{ }
\end{theorem}
\subsection{The $C^0$ version}
\begin{definition}
{\rm Let $f: (\mathbb{R}^n,0)\to (\mathbb{R}^p,0)$ be a $C^\infty$
map-germ. A $C^\infty$ deformation-germ $\Phi: (\mathbb{R}^n\times
\mathbb{R}^k,(0,0))\to (\mathbb{R}^p,0)$ of $f$ is said to be {\it
Thom trivial} if there exist $C$-regular stratifications in the
sense of \cite{bekka}, $\mathcal{S}$ of $\mathbb{R}^n\times
\mathbb{R}^k$, $\mathcal{T}$ of $\mathbb{R}^p\times \mathbb{R}^k$
and $\{\mathbb{R}^k\}$ of $\mathbb{R}^k$ such that the following two
conditions are satisfied:
\begin{enumerate}
\item The map-germ $(\Phi, \pi): (\mathbb{R}^n\times \mathbb{R}^k,(0,0))\to (\mathbb{R}^p\times \mathbb{R}^k,(0,0))$ is a Thom map-germ
with respect to $\mathcal{S}$ and $\mathcal{T}$, where $\pi:
(\mathbb{R}^n\times \mathbb{R}^k,(0,0))\to (\mathbb{R}^k,0)$ is the
canonical projection.
\item The canonical projection $\pi': (\mathbb{R}^p\times \mathbb{R}^k,(0,0))\to (\mathbb{R}^k,0)$ is a stratified map-germ with respect to
$\mathcal{T}$ and $\{\mathbb{R}^k\}$.
\end{enumerate}
}
\end{definition}
 It is known that there exists a controlled tube system for any $C$-regular stratification (see \cite{bekka}).
By Thom's second isotopy lemma, it follows that for any Thom trivial
deformation-germ $\Phi: (\mathbb{R}^n\times \mathbb{R}^k,(0,0))\to
(\mathbb{R}^p,0)$ of $f$, there exist germs of homeomorphisms $h:
(\mathbb{R}^n\times \mathbb{R}^k, (0,0))\to (\mathbb{R}^n\times
\mathbb{R}^k,(0,0))$ and $H: (\mathbb{R}^p\times \mathbb{R}^k,
(0,0))\to (\mathbb{R}^p\times \mathbb{R}^k,(0,0))$ such that the
diagram (1) in Section 1 is commutative.
\begin{theorem}[\cite{nishimuralondon}]\label{theorem 5}
Let $f, g: (\mathbb{R}^n,0)\to (\mathbb{R}^p,0)$ be two $C^\infty$
map-germs. Suppose that the following three conditions are
satisfied.
\begin{enumerate}
\item There exist a germ of $C^\infty$ diffeomorphism $s: (\mathbb{R}^n,0)\to (\mathbb{R}^n,0)$ and
a $C^\infty$ map-germ $M: (\mathbb{R}^n, 0)\to (GL(p, \mathbb{R}),
M(0))$ such that $f(x)=M(x)g(s(x))$.
\item The $C^\infty$ map-germ $F: (\mathbb{R}^n\times \mathbb{R}^p, (0,0))\to (\mathbb{R}^p,0)$ given by
\[
F(x, \lambda)=f(x)-M(x)\lambda
\]
is a Thom trivial deformation-germ of $f$.
\item The $C^\infty$ map-germ $G: (\mathbb{R}^n\times \mathbb{R}^p, (0,0))\to (\mathbb{R}^p,0)$ given by
\[
G(x, \lambda)=g(x)-M(s^{-1}(x))^{-1}\lambda
\]
is a Thom trivial deformation-germ of $g$.
\end{enumerate}
Then, $f$ and $g$ are topologically $\mathcal{A}$-equivalent.
\end{theorem}
A $C^\infty$ map-germ $f: (\mathbb{R}^n,0)\to (\mathbb{R}^p,0)$ is
said to be {\it MT-stable} if the jet extension of it is
multitransverse to the Thom-Mather canonical stratification of the
jet space (for the Thom-Mather canonical stratification, see for
instance \cite{gibsonetall, matherhowto, duplessiswall}). By
definition, it follows that any $C^\infty$ deformation of an
MT-stable map-germ is Thom trivial. Hence, the following
Fukuda-Fukuda theorem follows from Theorem \ref{theorem 5}.
\begin{corollary}[Fukuda-Fukuda theorem \cite{fukudafukuda}]
Let $f, g: (\mathbb{R}^n,0)\to (\mathbb{R}^p,0)$ be MT-stable
map-germs. Suppose that there exist a germ of $C^\infty$
diffeomorphism $s: (\mathbb{R}^n,0)\to (\mathbb{R}^n,0)$ and a
$C^\infty$ map-germ $M: (\mathbb{R}^n, 0)\to (GL(p, \mathbb{R}),
M(0))$ such that $f(x)=M(x)g(s(x))$. Then, $f$ and $g$ are
topologically $\mathcal{A}$-equivalent.
\end{corollary}
\begin{theorem}[\cite{nishimuralondon}]
Let $f, g: (\mathbb{R}^n,0)\to (\mathbb{R}^p,0)$ be two $C^\infty$
map-germs of rank zero. Suppose that the following two conditions
are satisfied.
\begin{enumerate}
\item There exist a germ of $C^\infty$ diffeomorphism $s: (\mathbb{R}^n,0)\to (\mathbb{R}^n,0)$ and
a $C^\infty$ map-germ $M: (\mathbb{R}^n, 0)\to (GL(p, \mathbb{R}),
M(0))$ such that $f(x)=M(x)g(s(x))$.
\item The $C^\infty$ map-germ $F: (\mathbb{R}^n\times \mathbb{R}^p, (0,0))\to (\mathbb{R}^p,0)$ given by
\[
F(x, \lambda)=f(x)-M(x)\lambda
\]
is a Thom trivial deformation-germ of $f$.
\end{enumerate}
Then, $f$ and $g$ are topologically $\mathcal{A}$-equivalent.
\end{theorem}
\section{Strategy}
\label{section 3}
The strategy to prove Theorems 1 and 2 is almost the same as the
strategies given in \cite{nishimuralondon, nishimuratopology}. In
\cite{nishimuralondon} (resp., \cite{nishimuratopology}), the
strategy is stated in terms of germs of homeomorphisms (resp., germs
of $C^r$ diffeomorphisms $(1\le r\le \infty)$).     It is possible
to be stated also in terms of germs of bi-Lipschitz homeomorphisms
as shown below.
\par
\medskip
Let $f,g: (\mathbb{R}^n,0)\to (\mathbb{R}^p,0)$ be two $C^\infty$
map-germs.  Suppose that there exist  a germ of $C^\infty$
diffeomorphism $s: (\mathbb{R}^n,0)\to (\mathbb{R}^n,0)$ and a
$C^\infty$ map-germ $M: (\mathbb{R}^n,0)\to (GL(p,\mathbb{R}),M(0))$
such that $f(x)=M(x)g(s(x))$. Then, we let $F:(\mathbb{R}^n\times
\mathbb{R}^p_\lambda, (0,0)\to (\mathbb{R}^p_y,0)$ be the map-germ
defined by
\[
F(x,\lambda)=f(x)-M(x)\lambda,
\]
where $\mathbb{R}^p_\lambda$ (resp., $\mathbb{R}^p_y$) is the
$p$-dimensional Euclidean space $\mathbb{R}^p$ which plays the role
of the parameter space (resp., the target space ) of $F$.    Suppose
furthermore that $F$ is Lipschitz trivial.  { }   Then, it
follows from Definition \ref{def1} that
there exist $L$-stratifications in the sense
of \cite{mostowski, parusinski, parusinski2}, $\mathcal{S}$ of
$\mathbb{R}^n\times \mathbb{R}^p_\lambda$, $\mathcal{T}$ of
$\mathbb{R}^p_y\times \mathbb{R}^p_\lambda$ and $\{\mathbb{R}^p_\lambda\}$ of
$\mathbb{R}^p_\lambda$ such that the three conditions in
Definition \ref{def1} are satisfied.
By the condition (iii) of Definition \ref{def1},
there { }exist
germs of bi-Lipschitz homeomorphisms $h: (\mathbb{R}^n\times
\mathbb{R}^p_\lambda, (0,0))\to (\mathbb{R}^n\times \mathbb{R}^p_\lambda,(0,0))$ and
$H: (\mathbb{R}^p_y\times \mathbb{R}^p_\lambda, (0,0))\to
(\mathbb{R}^p_y\times \mathbb{R}^p_\lambda,(0,0))$ such that the diagram (1) in Section 1 is
commutative. By the commutative diagram (1) we may put
\[
h(x, \lambda)=(h_1(x, \lambda), \lambda)\quad \mbox{and}\quad H(y,
\lambda)=(H_1(y,\lambda), \lambda).
\]
\begin{lemma}[\cite{nishimuratopology}, Lemma 2.1]\label{lemma 3.1}
$f(h_1(x, g(s(x))))=H_1(0,g(s(x))).$
\end{lemma}
\begin{lemma}[\cite{nishimuralondon}, Lemma 2.2]\label{lemma 3.2}
If the map germ from $(\mathbb{R}^p_\lambda,0)$ to
$(\mathbb{R}^p_y,0)$ defined by
\begin{equation}
\lambda\mapsto H_1(0,\lambda)
\end{equation}
is a germ of homeomorphism, then the map germ defined by
\begin{equation}
(x,\lambda)\mapsto (h_1(x,\lambda), H_1(0,\lambda))
\end{equation}
maps the germ of the set $(F^{-1}(0), (0,0))$ onto the germ of the
graph of $f$ at $(0,0)$.
\end{lemma}
\begin{lemma}
\label{lemma 3.3}
If the map germ (2)
is bi-Lipschitz,
then the endomorphism germ of
$(\mathbb{R}^n,0)$ defined by
\begin{equation}
x\mapsto h_1(x,g(s(x)))
\end{equation}
is also bi-Lipschitz.
\end{lemma}
{\it Proof of Lemma \ref{lemma 3.3}.}\qquad
The proof is almost the same as the proof of Lemma 2.3 in \cite{nishimuralondon}.
The map-germ (4) can be decomposed as follows.
\begin{equation}
x \mapsto  (x, g(s(x))) \mapsto (h_1(x,g(s(x))), H_1(0,g(s(x))))
\mapsto h_1(x,g(s(x))).
\end{equation}
If the map-germ (2) is bi-Lipschitz, then the map-germ (3) is also
bi-Lipschitz. Since the first map-germ of (5) is the germ of the
graph of $g\circ s$ at $(0,0)$, the composition of the first and the
second map-germ in (5) is also a germ of bi-Lipschitz homeomorphism
to the germ of the set
\begin{equation}
\left(\left\{
\left(h_1(x,g(s(x))), H_1(0,g(s(x)))\right)\; |\; x\in \mathbb{R}^n\right\}, (0,0)\right)
\end{equation}
with respect to the induced metric from $\mathbb{R}^n\times \mathbb{R}^p_y$.
The last map-germ in (5) is the restriction of the canonical projection
$(\mathbb{R}^n\times \mathbb{R}^p_y, (0,0))\to (\mathbb{R}^n,0)$ to (6).     By Lemma
\ref{lemma 3.2}, (6) is equal to the germ of the set
$(\mbox{\rm graph}(f), (0,0))$.
Hence, (4) is bi-Lipschitz.
\hfill\qquad $\Box$
\begin{lemma}\label{lemma 3.4}
Let $T_0$ be the stratum of $\mathcal{T}$ containing the origin $(0,0)$.
If $T_0$ is transverse to the linear space $\{0\}\times \mathbb{R}^p_\lambda$, then
the map-germ (2) is bi-Lipschitz.
\end{lemma}
{\it Proof of Lemma \ref{lemma 3.4}.}\qquad
It is clearly seen that the map-germ (2) is Lipschitz
even if $T_0$ is not transverse to $\{0\}\times \mathbb{R}^p_\lambda$.
\par
We show that if $T_0$ is transverse to $\{0\}\times \mathbb{R}^p_\lambda$, then
the map-germ (2) has its inverse which is also Lipschitz.
It is easily seen that the map-germ (2) has its inverse
if $T_0$ is transverse to $\{0\}\times \mathbb{R}^p_\lambda$.
For each $i$ $(1\le i\le p)$,
let $\eta_i+\frac{\partial}{\partial \lambda_i}$ be the germ of stratified Lipschitz
vector field such that the condition (iii) of Definition \ref{def1} is satisfied.
Moreover,
let $\Theta_i: (\mathbb{R}\times \mathbb{R}^p_y, (0,0))\to (\mathbb{R}^p_y, 0)$ be the germ of local flow for
$\eta_i$.        Then, the following holds
where $\lambda=(\lambda_1, \ldots, \lambda_p)\in \mathbb{R}^p_\lambda$.
\[
\Theta_1(\lambda_1; \Theta_2(\lambda_2; \ldots,
\Theta_p(\lambda_p, H_1(0, \lambda))\cdots )=0.
\]
For any sufficiently small $\lambda=(\lambda_1, \ldots, \lambda_p)$
in $\mathbb{R}^p_\lambda$ and $y\in \mathbb{R}^p_y$,
set $\Theta(\lambda, y)=\Theta_1(\lambda_1; \Theta_2(\lambda_2; \ldots,
\Theta_p(\lambda_p, y)\cdots )$.
Let $\lambda=(\lambda_1, \ldots, \lambda_p),
\widetilde{\lambda}=(\widetilde{\lambda}_1, \ldots, \widetilde{\lambda}_p)$ be sufficiently small two points of
$\mathbb{R}^p_\lambda$.     Then, since $\eta_i$ is Lipschitz for any $i$ $(1\le i\le p)$, there exists a positive constant $L$ such that the following holds for any $j$ $(1\le j\le p)$, where
$y_{j, \lambda}$ (resp., $y_{j, \widetilde{\lambda}}$) stands for
$ \Theta ((0, \ldots, 0, \lambda_j, \ldots, \lambda_p), H_1(0, \lambda))$
(resp., $\Theta ((0, \ldots, 0, \lambda_j, \ldots, \lambda_p), H_1(0, \widetilde{\lambda}))$) and
$y_{p+1, \lambda}$ (resp., $y_{p+1, \widetilde{\lambda}}$) stands for
$ H_1(0, \lambda))$ (resp., $H_1(0, \widetilde{\lambda}))$).
\begin{eqnarray*}
||
y_{j, \lambda}-y_{j, \widetilde{\lambda}}
||
& = &
||\Theta_j(\lambda_j, y_{{j+1}, \lambda})- \Theta_j(\lambda_j, y_{{j+1}, \widetilde{\lambda}})||
\\
{ } & \le &
||
y_{{j+1}, \lambda}-y_{{j+1}, \widetilde{\lambda}}
||
+
L
\left|
\int_{0}^{\lambda_j}
\left|\left|\Theta_{j}
(s, y_{j+1,\lambda}
)
-
\Theta_{j}
(s, y_{j+1,\widetilde{\lambda}})\right|\right|ds
\right|.
\end{eqnarray*}
Since $\lambda=(\lambda_1, \ldots, \lambda_p)$ is sufficiently small, we may assume that
$|\lambda_j|<1$ for any $j$ $(1\le j\le p)$.
Thus, by Gronwall's inequality (for instance, see \cite{hartman}), the following holds:
\begin{eqnarray*}\label{equation 1}
||\Theta(\lambda, H_1(0, \widetilde{\lambda})) ||
& = &
||\Theta(\lambda,
H_1(0, \lambda )) -
\Theta(\lambda, H_1(0, \widetilde{\lambda}))
|| \\
{ } & = &
||y_{1, \lambda}-y_{1, \widetilde{\lambda}} || \\
{ } & \le &
e||y_{2, \lambda}-y_{2, \widetilde{\lambda}} || \\
{ } & { }  &
\vdots \\
{ } & \le &
e^p||y_{{p+1}, \lambda}-y_{{p+1}, \widetilde{\lambda}} ||  \\
{ } & = &
 e^p||H_1(0, \lambda ) -
H_1(0, \widetilde{\lambda})
||.
\end{eqnarray*}
\par
For any $j$ $(1\le j\le p)$, set $\eta_j(y,
\lambda)=\sum_{i=1}^p\eta_{ij}(y, \lambda)\frac{\partial}{\partial
y_i}$, where $\frac{\partial}{\partial y_i}$ is the trivial vector
field of $\mathbb{R}^p_y$ with respect to the $i$-th coordinate
function $y_i$ of the standard coordinate neighborhood
$(\mathbb{R}^p_y, (y_1, \ldots, y_p))$. Let the following $p$ by $p$
matrix be denoted by $M(y,\lambda)$:
\[
M(y,\lambda)=
\left(
\begin{array}{ccc}
\eta_{11}(y, \lambda) & \cdots & \eta_{1p}(y, \lambda) \\
\vdots & \vdots & \vdots \\
\eta_{p1}(y, \lambda) & \cdots & \eta_{pp}(y, \lambda) \\
\end{array}
\right).
\]
Since $T_0$ is transverse to $\{0\}\times \mathbb{R}^p_\lambda$, we have that
$\eta_1(0,0), \ldots, \eta_p(0,0)$ are linearly independent.     This implies that
$\min_{||\widetilde{\lambda}||=1}||M(0,0)\widetilde{\lambda}||$ is positive.
By the continuity of
$\eta_i$ with respect to $y, \lambda\in \mathbb{R}^p$, it follows that
there exists a positive constant $M$ such that
\begin{equation}
0 < M\le \min_{||\widetilde{\lambda}||=1}||M(y,\lambda)\widetilde{\lambda}||,
\end{equation}
for any sufficiently small $y, \lambda\in \mathbb{R}^p$.
\par
Let $\varepsilon >0$ be a sufficiently small real number and let
$\psi=({}_1\psi, {}_2\psi): (-\varepsilon, \varepsilon)\to
\mathbb{R}^p_\lambda\times\mathbb{R}^p_\lambda$ be a regular curve
such that $\psi(0)\in \Delta$ and ${}_1\psi'(0)\ne {}_2\psi'(0)$,
where $\Delta$ is the diagonal set of
$\mathbb{R}^p_\lambda\times\mathbb{R}^p_\lambda$. Set $y_j\circ
{}_i\psi={}_i\psi_j$ $(1\le i\le 2, 1\le j\le p)$. Then, the
following equality holds:
\[
\Theta_1({}_2\psi_1(t) -{}_1\psi_1(t); \ldots,
\Theta_p({}_2\psi_p(t)-{}_1\psi_p(t);
\Theta({}_1\psi(t),
H_1 ( 0,{}_2\psi(t)))))\cdots )=0.
\]
By differentiating both sides of the above equality with respect to $t$ at $t=0$, we have
\[
M(0, {}_1\psi(0))\left({}_2\psi'(0)-{}_1\psi'(0)\right) +
\lim_{t\to 0}\frac{\Theta({}_1\psi(t),
H_1(0, {}_2\psi(t)))}{t}=0.
\]
Combining this equality and the Maclaurin expansion of ${}_2\psi(t)-{}_1\psi(t)$, we have the following:
\begin{equation}
\lim_{t\to 0}\frac{||\Theta({}_1\psi(t), H_1(0, {}_2\psi(t)))||}
{||{}_2\psi(t)-{}_1\psi(t)||}
=\left|\left|M(0, {}_1\psi(0))\frac{{}_2\psi'(0)-{}_1\psi'(0)}{||{}_2\psi'(0)-{}_1\psi'(0)||}\right|\right|.
\end{equation}
Let $\Psi: \mathbb{R}^p_\lambda\times \mathbb{R}^p_\lambda - \Delta\to \mathbb{R}$ be the function
defined by
\[
\Psi(\lambda, \widetilde{\lambda})=
\frac{||\Theta(\lambda, H_1(0, \widetilde{\lambda}))||}{||\widetilde{\lambda}-\lambda||}.
\]
Moreover, let $\widetilde{\pi}: B\to \mathbb{R}^p_\lambda\times \mathbb{R}^p_\lambda$ be the blowup
centered at $\Delta$.     The equality (8) shows that there exists the unique continuous function
$\widetilde{\Psi}: B\to \mathbb{R}$ such that $\widetilde{\Psi}=\Psi\circ \widetilde{\pi}$ on
$\mathbb{R}^p_\lambda\times \mathbb{R}^p_\lambda-\Delta$.
Thus, by the inequality (7), for any sufficiently small $\lambda, \widetilde{\lambda}\in \mathbb{R}^p_\lambda$ such that
$\lambda\ne \widetilde{\lambda}$, the following holds:
\[
\frac{M}{2}\le \frac{||\Theta(\lambda, H_1(0, \widetilde{\lambda}))||}
{||\widetilde{\lambda}-\lambda||}.
\]
\par
Therefore, for any sufficiently small $\lambda,
\widetilde{\lambda}\in \mathbb{R}^p_\lambda$, we have the following
inequality which shows that the inverse map-germ of (2) is actually
Lipschitz.
\[
||\lambda-\widetilde{\lambda}||\le e^p\frac{2}{M}||H_1(0, \lambda )
- H_1(0, \widetilde{\lambda})||.
\]
\hfill\qquad $\Box$

\section{Proofs of Theorem 1 and Theorem 2}
\label{section 4}
The proof of Theorem 1 (resp., Theorem 2) is almost the same as the
proof of Theorem 1.2 (resp., Theorem 1.3) in \cite{nishimuralondon}.
Theorem 1.3 in \cite{nishimuralondon} is proved by using
terms of Thom trivial deformations, $C$-regular stratifications,
Lemma 3.1 and the following two lemmas.
\begin{lemma}[\cite{nishimuralondon}, Lemma 2.3]
\label{lemma 4.1}
If the map germ from $(\mathbb{R}^p_\lambda,0)\to
(\mathbb{R}^p_y,0)$ defined by
\[
\lambda\mapsto H_1(0,\lambda)
\]
is a germ of homeomorphism,
then the endomorphism germ of
$(\mathbb{R}^n,0)$ defined by
\begin{equation}
x\mapsto h_1(x,g(s(x)))
\end{equation}
is also a germ of homeomorphism.
\end{lemma}
\begin{lemma}\label{lemma 4.2}
Suppose that $F$ is Thom trivial.  { }
Let $\mathcal{S}$, $\mathcal{T}$  and $\mathbb{R}^p_\lambda$ be the
$C$-regular stratifications of $\mathbb{R}^n\times \mathbb{R}^p_\lambda$,
$\mathbb{R}^p\times \mathbb{R}^p_\lambda$ and $\{\mathbb{R}^p_\lambda\}$ respectively
whose existence are guaranteed by Thom triviality of $F$.
Let $T_0$ be the stratum of $\mathcal{T}$ containing the origin $(0,0)$.
If $T_0$ is transverse to the linear space $\{0\}\times \mathbb{R}^p_\lambda$, then
the map-germ $\lambda\mapsto H_1(0,\lambda)$ is a germ of homeomorphism.
\end{lemma}
For the proof of Theorem 1.2 in \cite{nishimuralondon}, since we need to
consider perturbations of $F(x, \lambda)=f(x)-M(x)\lambda$ and
$g(x,\lambda)=g(x)-M(s^{-1}(x))^{-1}\lambda$ by matrices, two modifications of Lemma
\ref{lemma 4.2} are used.       By the proof of Lemma \ref{lemma 3.4},
it is not hard to prove that bi-Lipschitz version of these two modifications of
Lemma \ref{lemma 4.2} follows.
\par
Therefore, by replacing \lq\lq Thom trivial\rq\rq, \lq\lq
$C$-regular stratifications\rq\rq and Lemmas \ref{lemma 4.1} and
\ref{lemma 4.2} (resp., Lemma \ref{lemma 4.1} and two modifications
of Lemma \ref{lemma 4.2}) with \lq\lq Lipschitz trivial\rq\rq,
\lq\lq $L$-regular stratifications\rq\rq and Lemmas \ref{lemma 3.3}
and \ref{lemma 3.4} (resp., Lemma \ref{lemma 3.3} and the
bi-Lipschitz version of two modifications of Lemma \ref{lemma 4.2})
respectively, the Proof of  Theorem 1.3 (resp., Theorem 1.2) in
\cite{nishimuralondon} works well as the proof of Theorem 2 (resp.,
Theorem 1). \hfill\qquad $\Box$

\vspace{.5cm}

{ } \emph{\textbf{Remark 2.}} To give sufficient conditions
for bi-Lipschitz $\mathcal{A}$-equivalence of Lipschitz stable
map-germs is a key step in the attempt to answer the following  difficult open question in the subject of density of stable mappings:

\vspace{.2cm}

{\bf Question}: Is the set of Lipschitz stable mappings dense in $C^{\infty}(M^n,N^p),$ with the
Whitney $C^{\infty}$ topology, when
$(n,p)$ are outside the {\em nice dimensions} ?

\vspace{0.1cm}
It is well known that $C^0$ stable maps are dense in general and that
$C^1$ stability fails to be dense in the function space exactly when
$C^\infty$ stability fails to be dense, that is, outside the nice
dimensions. However, few results are known about Lipschitz stable maps. Actually we know of no example of a Lipschitz stable map
(or germ) which is not $C^\infty$ stable.

 { }

\subsection*{Acknowledgement}
The authors are grateful to the reviewers for careful reading of the first draft of their paper and giving
important comments.

{{ }
\begin{flushright}



\end{flushright}
}
\end{document}